\documentclass[12pt]{article}
\usepackage{latexsym,amssymb}
\font\msxm=msxm10 scaled 1200
\font\of=msbm10 scaled 1200
\def\lesssim{\mbox{\hspace{1mm}\msxm\char'56}\hspace{1mm}}
\def\R{\mbox{\of R}}
\def\C{\mbox{\of C}}
\def\Z{\mbox{\of Z}}
\def\N{\mbox{\of N}}

\begin{document}
\hoffset =-1cm
\voffset =-2cm

\begin{center}
\noindent
{\Large\bf Two dimensional Fuchsian systems

\vspace{1ex}
and the Chebyshev property
\footnote {Research partially supported
by a grant from the NSF of Bulgaria
and CNRS France.}}

\vspace{1cm}
\noindent
{\bf Lubomir Gavrilov}

\vspace{1ex}
Laboratoire de Mathematiques Emile Picard,  Univ. Paul Sabatier

\noindent
31062 Toulouse Cedex 04, France, e-mail: gavrilov@picard.ups-tlse.fr
\vspace{1cm}
\noindent

{\bf Iliya  D. Iliev}

\vspace{1ex}
Institute of Mathematics, Bulgarian Academy of Sciences

\noindent
P.O. Box 373, 1090 Sofia, Bulgaria, e-mail: iliya@math.bas.bg

\vspace{5mm}
\noindent
October 26, 2001
\end{center}

\vspace{5mm}
\noindent
{\bf Abstract.}
Let $(x(t),y(t))^\top$ be a solution of a Fuchsian system of order two with
three singular points. The vector space of functions of the form
$P(t)x(t)+Q(t)y(t)$, where $P,Q$ are real polynomials, has a natural filtration
of vector spaces, according to the asymptotic behaviour of the functions at
infinity. We describe a two-parameter class of Fuchsian systems, for which
the corresponding vector spaces obey  the Chebyshev property (the maximal
number of isolated zeroes of each function is less than the dimension of the
vector space).

Up to now, only a few particular systems were known to possess such a
nonoscillation property. It is remarkable that most of these systems are of
the type studied in the present paper. We apply our
results in estimating the number of limit cycles that appear after small
polynomial perturbations of several quadratic or cubic Hamiltonian
systems in the plane.

\vspace{5mm}
\noindent
2000 MSC scheme numbers: 34C07, 34C08, 34C05

\vspace{10mm}
\normalsize
\noindent
{\large\bf 1. Introduction}

\vspace{2ex}
\noindent
In many bifurcation problems the main difficulty is to estimate the number
of isolated zeroes of certain functions of the form
$$I(h)=p_1(h)I_1(h)+p_2(h)I_2(h),\quad h\in\Sigma,\eqno(1)$$
where $p_1(h)$ and $p_2(h)$ are polynomials,
and the vector function ${\bf I}(h)=(I_1(h), I_2(h))^\top$
satisfies a two-dimensional first order Fuchsian system
$${\bf I}(h)={\bf A}(h){\bf I}'(h),\qquad '=d/dh, \eqno(2)$$
with a first degree polynomial matrix ${\bf A}(h)$. Typically,
$I_1(h)$ and $I_2(h)$ are complete Abelian integrals along the ovals
$\delta(h)$ within a continuous (in $h$) family of ovals contained in the
level sets of a fixed real polynomial $H(x,y)$ (called the Hamiltonian),
and $\Sigma\subset\R$ is the maximal open interval of existence of
such ovals $\delta(h)$. See Table 1.

In the present paper, our main assumptions on (2) are the following:

\vspace{2ex}
\noindent
(H1) {\it ${\bf A}'$ is a constant matrix having real distinct
eigenvalues.}

\vspace{1ex}
\noindent
(H2) {\it The equation ${\rm det}\,{\bf A}(h)=0$ has real distinct roots
$h_0$, $h_1$ and the identity
\hspace*{1cm}${\rm trace}\,{\bf A}(h)\equiv({\rm det}\,{\bf A}(h))'$
holds.}

\vspace{1ex}
\noindent
(H3){\it $\;{\bf I}(h)$ is analytic in a neighborhood of $h_0$.}

\vspace{2ex}
\noindent

The conditions that ${\bf A}'$ is a constant matrix and
${\rm det}\,{\bf A}(h)$ has distinct roots imply that the singular
points of the system
$${\bf I}'(h) = {\bf A}^{-1}(h){\bf I}(h)$$
(including $\infty$) are regular, i.e. it is of Fuchs type. Further, the
condition ${\rm trace}\,{\bf A}(h)\equiv({\rm det}\,{\bf A}(h))'$
implies that
the characteristic exponents of (2) at $h_0$ and $h_1$ are $\{0,1\}$.
In the formulation of our main result below, we assume for definiteness
that $h_0<h_1$. A similar result holds if $h_0>h_1$.
Clearly if $h_0<h_1$, and the function ${\bf I}(h)$ is analytic in a
neighborhood of $h=h_0$, then it also possesses an analytic
continuation in the complex domain $\C \backslash [h_1,\infty)$.

\vspace{2ex}
\noindent
{\bf Definition 1.} {\it
The real vector space of functions $V$ is said to be Chebyshev
in the complex domain ${\cal D}\subset \C$ provided
that every function $I\in V\setminus\{0\}$ has at most ${\rm dim}\,V-1$
zeros in $\cal D$.
$V$ is said to be Chebyshev with accuracy $k$ in $\cal D$ if any
function $I\in V\setminus\{0\}$ has at most $k+{\rm dim}\, V -1$ zeros
in $\cal D$. }

\vspace{2ex}
\noindent
{\bf Definition 2.} {\it Let $I(h)$, $h\in \C$ be a function, locally
analytic in a neighborhood of $\infty$, and $s\in\R$. We shall write
$I(h) \lesssim h^s$, provided that for every sector $S$ centered at
$\infty$ there exists a non-zero constant $C_S$ such that
 $|I(h)| \leq C_S |h|^s$ for all sufficiently big $|h|$, $h\in S$.}

\vspace{2ex}
\noindent
For systems (2) satisfying (H1) and (H2), the
characteristic exponents at infinity are
$-\lambda$ and $-\mu$ where $\lambda'=1/\lambda$ and $\mu'=1/\mu$ are the
eigenvalues of the constant matrix ${\bf A}'$. According to (H2),
$\lambda+\mu=2$. Let us denote $\lambda^*=2$ if $\lambda$ is integer and
$\lambda^*={\rm max}\,(|\lambda-1|, 1-|\lambda-1|)$ otherwise.

Take $s\geq \lambda^*$ and consider the real vector space of functions
 $$V_s=\{I(h)= P(h)I_1(h)+Q(h)I_2(h):\; P,Q\in \R[h],\; I(h) \lesssim
h^s \}$$
where ${\bf I}=(I_1(h),I_2(h))^\top$ is a non-trivial  solution of (2),
holomorphic in a neighborhood of $h=h_0$. As $\lambda,\mu\not\in\{0,1,2\}$,
 the vector function
${\bf I}(h)$ is uniquely determined, up to multiplication by a constant,
and $I_1(h_0)=I_2(h_0)=0$ (see Proposition 1). Clearly, $V_s$ is invariant
under linear transformations in (2) and affine changes of the argument $h$.
The restriction $s\geq \lambda^*$ is taken to guarantee that $V_s$ is not
empty.

Recall that $h_0<h_1$ are the roots of ${\rm det}\,A(h)=0$.
 Our main result in section 2 is the following.

\vspace{2ex}
\noindent
{\bf Theorem 1.} {\it
 Assume that conditions {\rm (H1) -- (H3)} hold.
 If $\lambda\not \in \Z$,
 then $V_s$ is a Chebyshev vector space with accuracy
 $1+[\lambda^*]$ in the complex domain
 ${\cal D}=\C \backslash [h_1,\infty)$. If $\lambda \in
\Z$, then $V_s$ coincides with the space of real
polynomials of degree at most $[s]$ which vanish at $h_0$ and $h_1$.}

\begin{center}
\begin{tabular}{|c|c|c|c|c|}\hline
No. &  $H$, $\;\bf I$  & $\Sigma$  &  $\bf A$  &  ${\rm det}\,{\bf A},\;
  {\rm tr}\,{\bf A}$ \\ \hline
    &      &   & &\\[-5mm]
1&$\begin{array}{c} H=y^2+x^2-x^3  \\
{\bf I}=(\int_{H=h}ydx, \int_{H=h}xydx)\end{array}$ &
 $(0,\frac{4}{27})$ &
 $\left(\!\!\begin{array}{cc}\frac65h& -\frac{4}{15}\\ \frac{4}{35}h
 &\frac67h-\frac{16}{105}\end{array}\!\!\right)$  &
$\begin{array}{l}\frac{36}{35}h^2-\frac{16}{105}h\\
\frac{72}{35}h-\frac{16}{105}\end{array}$ \\
    &      &   & & \\[-5mm]
  \hline
    &      &   & & \\[-5mm]
2&$\begin{array}{c} H=y^2+x^2-xy^2\\{\bf I}=(\int_{H=h}ydx, \int_{H=h}xydx)
\end{array}$ &
 $(0,1)$ &
 $\left(\!\!\begin{array}{cc}\frac43h&-\frac43\\ \frac{4}{15}h&
 \frac45h-\frac{16}{15}\end{array}\!\!\right)$  &
$\begin{array}{l}\frac{16}{15}h^2-\frac{16}{15}h\\
\frac{32}{15}h-\frac{16}{15}\end{array}$\\[5mm]
    &      &   & & \\[-5mm]
  \hline
    &      &   & & \\[-5mm]
3&$\begin{array}{c} H=\frac12y^2+\frac12x^2-\frac13x^3+xy^2\\
{\bf I}=(\int_{H=h}ydx,
\int_{H=h}x^2ydx)\end{array}$ &
 $(0,\frac16)$ &
 $\left(\!\!\begin{array}{cc}\frac32h&-\frac12\\ \frac{3}{16}h&
 \frac34h-\frac{3}{16}\end{array}\!\!\right)$  &
$\begin{array}{l}\frac{9}{8}h^2-\frac{3}{16}h\\
\frac{9}{4}h-\frac{3}{16}\end{array}$\\
    &      &   & & \\[-5mm]
  \hline
    &      &   & & \\[-5mm]
4&$\begin{array}{c} H=y^2+x^2+x^4\\{\bf I}=(\int_{H=h}ydx, \int_{H=h}x^2ydx)
\end{array}$ &
 $(0,\infty)$ &
 $\left(\!\!\begin{array}{cc}\frac43h&-\frac23\\-\frac{2}{15}h&
 \frac45h+\frac{4}{15}\end{array}\!\!\right)$  &
$\begin{array}{l}\frac{16}{15}h^2+\frac{4}{15}h\\
\frac{32}{15}h+\frac{4}{15}\end{array}$\\
    &      &   & & \\[-5mm]
  \hline
    &      &   & & \\[-5mm]
5&$\begin{array}{c} H=y^2+x^2-x^4\\{\bf I}=(\int_{H=h}ydx, \int_{H=h}x^2ydx)
\end{array}$ &
 $(0,\frac14)$ &
 $\left(\!\!\begin{array}{cc}\frac43h&-\frac23\\ \frac{2}{15}h&
 \frac45h-\frac{4}{15}\end{array}\!\!\right)$  &
$\begin{array}{l}\frac{16}{15}h^2-\frac{4}{15}h\\
\frac{32}{15}h-\frac{4}{15}\end{array}$\\
    &      &   & & \\[-5mm]
  \hline
    &      &   & & \\[-5mm]
6&$\begin{array}{c} H=y^2+x^2+x^2y^2\\{\bf I}=(\int_{H=h}ydx, \int_{H=h}x^2ydx)
\end{array}$ &
 $(0,\infty)$ &
 $\left(\!\!\begin{array}{cc}2h&-2\\ -\frac{2}{3}h&\frac23h+\frac{4}{3}
 \end{array}\!\!\right)$  &
$\begin{array}{l}\frac43h^2+\frac43h\\ \frac83h+\frac43\end{array}$\\
    &      &   & & \\[-5mm]
  \hline
    &      &   & & \\[-5mm]
7&$\begin{array}{c} H=y^2+x^2-x^2y^2\\{\bf I}=(\int_{H=h}ydx, \int_{H=h}x^2ydx)
\end{array}$ &
 $(0,1)$ &
 $\left(\!\!\begin{array}{cc}2h&-2\\ \frac{2}{3}h&\frac23h-\frac{4}{3}
 \end{array}\!\!\right)$  &
$\begin{array}{l}\frac43h^2-\frac43h\\ \frac83h-\frac43\end{array}$\\
   &      &   & & \\[-5mm]
   \hline
 &      &   &  &\\[-5mm]
8& $\begin{array}{c}
H=x^{-3}(y^2-2x^2+x)\\
{\bf I}=(\int_{H=h}x^{-3}ydx, \int_{H=h}x^{-4}ydx)
\end{array}$ &
 $(-1,0)$ &
 $\left(\!\!\begin{array}{cc}\frac43h & \frac43\\ \frac{4}{15}h
 &\frac45h+\frac{16}{15}
 \end{array}\!\!\right)$&
 $\begin{array}{l}\frac{16}{15}h^2+\frac{16}{15}h\\ \frac{32}{15}h+
 \frac{16}{15}\end{array}$ \\[5mm]
 \hline
  \end{tabular}

\vspace{5mm}
{\small Table 1. Examples of systems for integrals ${\bf I}=(I_1, I_2)$ and

Hamiltonian functions $H$ which satisfy hypotheses (H1) -- (H3).}
\end{center}

As an application of Theorem 1 let us consider a polynomial
perturbation of a planar Hamiltonian system
$$\begin{array}{l}
\dot{x}=H_y+\varepsilon f(x,y),\\
\dot{y}=-H_x+\varepsilon g(x,y),\end{array}\eqno(3)$$
where $\varepsilon$ is a small parameter, the degree of the polynomials
$f,g$ does not exceed $n$ and $H$ is some of the Hamiltonians
from Table 1. Define the function
$$h\;\to\;I(h)=\oint_{H=h}[g(x,y)dx-f(x,y)dy],\quad
h\in\Sigma.\eqno(4)$$
As is well known, if $I(h)\not\equiv 0$ in $\Sigma$, then the number of
limit cycles in (3) bifurcating for small $\varepsilon$ from the
periodic
orbits of the unperturbed Hamiltonian system is bounded by the number of
isolated zeroes of $I(h)$ in $\Sigma$. Define the linear space
${\cal V}_n$ of integrals given by (4) for ${\rm deg}\,f,g\leq n$.
Denote by $h_1$ the nonzero critical value of the Hamiltonian and by
${\cal D}$ the complex plane cut along the part of the real axis between
$h_1$ and $\infty$ not containing the other critical value $h_0=0$.
Then applying Theorem 1, we obtain the following results.

\vspace{2ex}
\noindent
{\bf Theorem 2.} {\it For each of the systems $1)-5)$ in Table 1,
 the linear space of integrals ${\cal V}_n$ is Chebyshev with accuracy
 one in ${\cal D}$. In particular, ${\cal V}_n$ is Chebyshev in $\Sigma$.}

\vspace{2ex}
\noindent
{\bf Theorem 3.} {\it For systems $6)$ and $7)$ in Table 1, the linear
space of integrals ${\cal V}_n$ is Chebyshev with accuracy
one in ${\cal D}$, if $n\leq 6$, and with accuracy $[\frac{n+1}{4}]$,
if $n\geq 7$. In particular, ${\cal V}_n$ is Chebyshev in
$\Sigma$, if $n\leq 6$,  and Chebyshev with accuracy $[\frac{n-3}{4}]$,
if $n\geq 7$.}

\vspace{2ex}
\noindent
Roughly speaking, Theorems 2 and 3 imply that, for the systems 1)--7)
from Table 1, the number of limit cycles in (3) born out of periodic
orbits under small polynomial perturbations which are transversal to
the integrable directions, is less than the dimension of the linear
space of these perturbations (with certain accuracy if $n\geq 7$ in
cases 6) and 7)). Clearly, a bound obtained by establishing the Chebyshev
property, is always the optimal one.

Case 8) from Table 1 is non-Hamiltonian one and requires slightly
different approach. See the end of the paper for results about it.

Let us recall that Theorem 2 in case 1) was proved earlier by Petrov [10].
Some less general (or a little bit different) results concerning
cases 3) -- 5) can be found in [1], [5], [6] and [8].

\vspace{3ex}
\noindent
{\large\bf 2. The Chebyshev property}

\vspace{2ex}
\noindent
We intend first to obtain a normal form for the matrices satisfying (H1)
and
(H2). For this purpose, we perform in (2) a linear transformation
bringing
${\bf A}'$ to a diagonal form and then translate the critical value
$h_0$ to
the origin. The matrix in (2) takes the form
$${\bf A}(h)=\left(\begin{array}{cc}
{\displaystyle\frac{2h-h_1}{2\lambda}}&
{\displaystyle\frac{\omega h_1}{2\lambda}}\\[3mm]
{\displaystyle\frac{h_1}{2\mu\omega}} &
{\displaystyle\frac{2h-h_1}{2\mu}}
\end{array}\right)\eqno(5)$$
where $h_1$ is the nonzero critical value and $\omega$ is a free
parameter.
This is the normal form we will use in this section.
In applications, another normal form takes place. To obtain it,
we apply additional linear transformation in (5)
$(I_1,I_2)\rightarrow (I_1, I_1/\omega+I_2)$ bringing ${\bf A}(h)$ to
$${\bf A}(h)=\left(\begin{array}{cc}
{\displaystyle\frac{h}{\lambda}}&
{\displaystyle\frac{\omega h_1}{2\lambda}}\\[3mm]
{\displaystyle\frac{(\lambda-\mu)h}{\lambda\mu\omega}} &
{\displaystyle\frac{h}{\mu}-\frac{h_1}{\lambda\mu}}
\end{array}\right).\eqno(6)$$
Evidently, equations (5) and (6) present three-parameter
families of matrices which can be reduced to two-parameter ones by
moving
$h_1$ to $1$. We note that all the examples in Table 1
are taken in the normal form (6), with $\frac12\leq\lambda<\mu\leq\frac32$.

Before to prove Theorem 1, we need some preparation. Without any loss
of generality we may use the normal form (5), with $h_1=1$. Hence, we will
consider $t=(h-h_0)/(h_1-h_0)$ as the argument and will assume throughout
this section that (2) is rewritten as a system
${\bf I}(t)={\bf A}(t){\bf I}'(t)$ for
${\bf I}(t)=(x(t),y(t))^\top\equiv(I_1(h),I_2(h))^\top$, with
$${\bf A}(t)=\left(\begin{array}{cc}
{\displaystyle\frac{2t-1}{2\lambda}}&
{\displaystyle\frac{\omega}{2\lambda}}\\[3mm]
{\displaystyle\frac{1}{2\mu\omega}}&
{\displaystyle\frac{2t-1}{2\mu}}\end{array}\right).\eqno(7)$$

\vspace{2ex}
\noindent
{\bf Proposition 1.} {\it The functions $x(t)=I_1(h)$ and $y(t)=I_2(h)$
satisfy
equations}
$$t(t-1) x'' = \lambda(\lambda-1) x, \eqno{(8)}$$
$$ t(t-1) y'' = \mu(\mu-1) y. \eqno{(9)}$$
{\bf Proof.} The most easy proof is a straightforward calculation
which we left to the reader (cf. [8]).

\vspace{2ex}
\noindent
{\bf Proposition 2.} {\it Let $\lambda\neq 0,1$ and $x(t)$ be a
nontrivial solution of $(8)$ which is analytic in a neighborhood of $t=0$
$($or $t=1)$. Then $x(t)\neq 0$ for $t<0$ $($respectively,
for $t>1)$. In particular, if $\lambda\in\Z$,
then $x(t)$ is a special kind of ultra-spherical polynomial and has
all of its zeros in the interval $[0,1]$.}

\vspace{2ex}
\noindent
{\bf Proof.} The assertion is well known for $\lambda$ integer. In this
case $x(t)$ is a kind of ultra-spherical (Gegenbauer) polynomial
$[12]$ of degree $\lambda$  if $\lambda\geq 2$ and of degree $1-\lambda$
if $\lambda\leq -1$. Although the result might be known for $\lambda$
not integer too, we will for completeness give the proof for this case.
Let $x(t)$ be analytic near $t=0$ (the other case is similar).
Take the function
$$z(t)=\frac{t^2-t}{2-\lambda}x'+\frac{1-\lambda t}{2-\lambda}x.$$
Then $z'=tx'-\lambda x$ and $x'$, $z'$ together satisfy a system
$$\begin{array}{l}
(t^2-t)x''=(\lambda-1)(tx'-z')\\
(t^2-t)z''=(\lambda-1)(tx'-tz').\end{array}$$
As $x'(0)\neq 0$ and $z'(0)=0$, the ratio $w=z'/x'$
is an analytical function in a neighborhood of $t=0$
satisfying the Riccati equation
$$\frac{t^2-t}{\lambda-1}w'(t)=w^2(t)-2tw(t)+t$$
and $w(0)=0$.
Consider in the $(t,w)$-plane the zero isocline given by the hyperbola
$w^2-2tw+t=0.$ It goes through the origin and has a vertical asymptote
at that point. It is easy to conclude that for $t<0$, the graphic of $w$
is placed inside the left branch of the hyperbola and either
$w(t)>0$ or $w(t)<0$ for all $t<0$, depending on whether $w'(0)$ is
negative or positive. Therefore $x'(t)$ and $z'(t)$ do not change signs
for $t<0$. As $x(0)=0$, the assertion follows.   $\Box$

\vspace{2ex}
\noindent
{\bf Proposition 3.} {\it Let $\lambda<1$ and $x(t)$ be a nontrivial
solution of {\rm (8)} which is  analytic in a neighborhood of $t=0$.
If $\lambda \not\in \Z$, then $x(t)$ has at most
$1+[\lambda^*]$ zeros in the complex domain
${D} = \C\backslash [1,\infty)$.}

\vspace{2ex}
\noindent
{\bf Proof.} Consider the analytic continuation of $x(t)$ in the complex
domain ${D}=\C\backslash [1,\infty)$. We shall count the zeros of
$x(t)$ in ${D}$ by making use of the argument principle.
Let $R$ be a big enough constant and $r$ a small enough constant.
 Denote by $\tilde{D}$ the set obtained  by removing
 the small disc $\{|t-1|<r\}$ from ${D}\cap\{|t|<R\}$.
 To estimate the number of the zeros of $x(t)$ in $\tilde{D}$,
 we shall evaluate the increment $\Delta_{\partial \tilde{{D}}}
{\rm Arg}\,x(t)$ of the argument of the function $x(t)$ along the
boundary
 of $\tilde{D}$, traversed in a positive direction. Then, according
 to the argument principle, we have that the number of the zeros of
$x(t)$ in $\tilde{D}$ equals
 $$\frac{\Delta_{\partial \tilde{D}} {\rm Arg}\, x(t)}{2\,\pi}.$$
 The monodromy group of the  equation in (8) is reducible if and only
 if $\lambda \in \Z$ [4, Theorem 4.3.2]. Therefore, if $\lambda \not\in
\Z$, then in a neighborhood of $t=1$ we have
$$x(t)= \xi (t) \log(t-1) +\eta (t)$$
where $\xi(t)$, $\eta(t)$ are analytic in a neighborhood of $t=1$,
$\xi(t)$ is a non-trivial solution of (8), $\xi (1)=0$. Moreover, a local
analysis shows that $ \lim_{t\rightarrow 1^-}x(t)=\eta(1)=const\neq 0$.
Therefore the increase of the argument of $x(t)$, when running the boundary
of $\{|t-1|<r\}$, is close to zero. Along the half line $(1,\infty)$ the
imaginary part of $x(t)$ equals $\pi\xi(t)$ which does not vanish,
by Proposition 2. Finally if $|t|$ is sufficiently
big then we have
$$|x(t)|\leq c|t|^\lambda\;\,{\rm if}\;\,\lambda>{\textstyle\frac12},\quad
|x(t)|\leq c|t|^{1-\lambda}\;\,{\rm if}\;\,\lambda<{\textstyle\frac12},\quad
|x(t)|\leq c|t^\frac12\log t|\;\,{\rm if}\;\,\lambda={\textstyle\frac12},$$
where $c$ is a non-zero constant. The increase of the argument of $x(t)$,
when running the boundary of $\{|t|<R\}$ is close to $2\pi \lambda^*.$
Summing up the above information, we obtain that the increase of
the argument of $x(t)$, when running the boundary of ${D}$, is at most
$2\pi+ 2\pi \lambda^* $.  We conclude that $x(t)$ has at most
$1+[\lambda^*]$ zeros in ${D}$ which completes the proof of Proposition 3.
$\Box$

\vspace{1ex}
We also need a more detailed information about the structure of the linear
space $V_s$ and an explicit formula for ${\rm dim}\, V_s$. The only
interesting case is when $\lambda$ and $\mu$ are not integer.

\vspace{2ex}
\noindent
{\bf Proposition 4.} {\it Let $s\geq \lambda^*$ and $\lambda,\mu$
be not integer. Then}
$${\rm dim}\;V_s=\left\{\begin{array}{l}
2s-1,\;{\rm if}\;\lambda-\mu \;{\rm and} \; s-\frac12\;\;
{\rm are}\;{\rm integer},\\[2mm]
[s-\lambda]+[s-\mu]+2,\;{\rm otherwise.}
\end{array}\right.$$

\vspace{2ex}
\noindent
{\bf Proof.} Without loss of generality, we can use the coordinates in
which $\bf A$ takes a form (7) and ${\bf I}(h)=(x(t),y(t))^\top$. To
reduce the number of cases,
let us assume that $\lambda>\mu$ (when $\lambda<\mu$, the analysis is
similar).

We begin our analysis with the case when $s\geq\lambda$.
Assume first that $\lambda-\mu$ is not integer. Then one can take any
solution of (2) near infinity in the form
$${\bf I}=\pmatrix{x\cr y}=
a\pmatrix{t^\lambda-\frac{\lambda}{2}t^{\lambda-1}+\ldots\cr
\alpha t^{\lambda-1} +\ldots }
+b\pmatrix{\beta t^{\mu-1} +\ldots \cr
t^\mu-\frac{\mu}{2}t^{\mu-1}+\ldots}$$
where
$$\alpha=\frac{\lambda}{2\omega(\mu-\lambda+1)},\;
\beta=\frac{\mu\omega}{2(\lambda-\mu+1)}.$$
Since ${\bf I}$ is analytic in a neighborhood of zero, the constants
$a$ and $b$ are both nonzero. Indeed, if $ab=0$ then ${\bf I}$ defines
an
 one-dimensional subspace in the space of all solutions, which is
invariant
under the monodromy group of (2), and hence of (8),(9). This is however
impossible, as the latter groups are irreducible for $\lambda ,\mu \not
\in \Z$.

Given $s\geq\lambda$, then the function $I(h)$ in the
definition of $V_s$ contains monomials of the form
$t^kx$, $0\leq k\leq K$, $t^ly$, $0\leq l\leq L$, where
$K\leq{\rm min}\,(s-\lambda,s-\mu+1)$, $L\leq{\rm
min}\,(s-\lambda+1,s-\mu)$.
Using that $\lambda+\mu=2$ and $\lambda>\mu$, one obtains
$K\leq s-\lambda+{\rm min}\,(0, 2\lambda-1)=s-\lambda.$ Similarly,
$L\leq s-\mu+{\rm min}\,(0,2\mu-1)=s-\lambda+1$ if $\lambda-\mu>1$ and
$L\leq s-\mu$ otherwise.

Among these monomials, other special combinations may be
involved in $V_s$ if $\lambda-\mu>1$. Define the functions
$z_1=ty-\alpha_1 x$, $z_m=tz_{m-1}-\alpha_m x$, $m\geq 2$,
where $\alpha_1=\alpha$ and the constant $\alpha_m$ is determined so
that
the coefficient at $t^\lambda$ in $z_m$ is zero. Denote $M=[s-\mu]-K-1$.
Clearly, then $t^{K+1}z_m\in V_s$ for $1\leq m\leq M$.
Moreover, any combination $t^{K+1}(P(t)x+Q(t)y)$ which belongs to $V_s$
is a linear combination of the ``monomials" $t^{K+1}z_m$.

Thus, ${\rm dim}\, V_s= K+L+2$ for $|\lambda-\mu|<1$ and
${\rm dim}\, V_s= K+L+M+2$ otherwise, which yields
${\rm dim}\, V_s= [s-\lambda]+[s-\mu]+2$ in both cases.

Assume now that $\lambda-\mu$ is integer but $\lambda$ and $\mu$ are
not.
If $\lambda-\mu>1$,
then one can take any solution of (2) near infinity in the form
$${\bf I}= a\pmatrix{t^\lambda-\frac{\lambda}{2}t^{\lambda-1}+\ldots\cr
\alpha t^{\lambda-1}+\ldots}
+(a\gamma\log t+b)\pmatrix{\beta t^{\mu-1}+\ldots\cr
t^\mu-\frac{\mu}{2}t^{\mu-1}+\ldots},\;\;\gamma\neq 0.$$
As in the previous case, this yields $K=[s-\lambda]$,
$L=[s-\lambda+1]=K+1$
and $M=[s-\mu]-K-1$ if $s-\frac12$ is not integer,  $M=s-\mu-K-2$ if
$s-\frac12$ is integer. In the first case we obtain the same result as
above,
and in the second case ${\rm dim}\,V_s=K+L+M+2=2s-1$.

Finally, if $\lambda=\frac32$, $\mu=\frac12$, we have respectively
$${\bf I}=a\pmatrix{t^\frac32-\frac34t^\frac12-\frac{9}{64}t^{-\frac12}
+\ldots\cr \frac{3}{8\omega}t^{-\frac12}+\ldots}
+(-\frac{3a}{4\omega}\log
t+b)\pmatrix{\frac{\omega}{8}t^{-\frac12}+\ldots\cr
t^\frac12-\frac14t^{-\frac12}+\ldots}.$$
Clearly $K=[s-\frac32]$, $L=[s-\frac12]$ if $s-\frac12$ is not integer
and $L=K$ otherwise.
Since no other combinations are involved in $V_s$ in this case, the
result follows immediately.

In the case when $\lambda>s\geq \lambda^*$,
the analysis is simpler. We use the same formulas for ${\bf I}$ as above.
One has either (a) $\lambda>s\geq\lambda-1$ and
$\lambda-\mu\geq 1$, or (b) $\lambda>s\geq \mu$ and $\lambda-\mu< 1$.
If $s=\mu=\frac12$, then $V_s$ is empty.
In all other cases, $y\in V_s$. In case (b), $V_s$ contains no other
functions. In case (a), if $s\geq \mu+1$ and $\lambda-\mu$ is not
integer, then also $z_m\in V_s$ for $1\leq m\leq [s-\mu]$.
The same is true if $\lambda-\mu$ is integer but $s-\frac12$ is not.
Finally, if both $\lambda-\mu$ and $s-\frac12$ are integer, and $s\geq \mu+2$,
then $V_s$ contains the functions  $z_m$, $1\leq m\leq [s-\mu]-1$.
Clearly, in all the cases above we obtain a formula for ${\rm dim}\,V_s$
as asserted. $\Box$

\vspace{2ex}
\noindent
{\bf Proof of Theorem 1.} For integer $\lambda,\mu$ the assertion
is obvious since $I_1$ and $I_2$ are different ultra-spherical polynomials
which have no common zeros except the simple ones at $h_0$ and $h_1$.

Assume below that $\lambda ,\mu \not \in \Z$ and let $\lambda>\mu$
(for definiteness). Suppose as before that the matrix ${\bf A}$ takes
the form (7), and let $I(t)= P(t)x (t)+Q(t)y(t)\in V_s$, where
$(x(t)=I_1(h)$, $y(t)=I_2(h))$ is the holomorphic solution
of (2) vanishing at the origin. When $P(t)\equiv 0$, the assertion is
evident. When $P\not\equiv 0$, we use again the argument principle
to count the zeros of $I(t)$ in $D=\C\backslash [1,\infty)$.
Consider in $D$ the meromorphic function
$$F(t)= P(t)\frac{x(t)}{y(t)}+ Q(t).$$
Below we calculate the increase of its argument when running the boundary
of $D$. The local structure of the solutions of (8), (9) in a
neighborhood of $t=1$ implies that $\lim_{t\rightarrow 1}x(t) \neq 0$,
$\lim_{t\rightarrow 1}y(t)\neq 0$. Therefore the increase
of the argument of $F(t)$, when running the boundary of $\{|t-1|< r\}$
is close to zero. As $x(t)$, $y(t)$ are real-analytic on $(-\infty,1)$,
then along the half-line $(1,\infty) $
$${\rm Im}\, F(t) = P(t)\,{\rm Im} \frac{x(t)}{y(t)}=P(t)
\frac{{\rm det}\left(\begin{array}{rl}\overline{y(t)} & y(t)\\
\overline{x(t)} & x(t) \end{array}\right)}{2i\,|y(t)|^2}.$$
As $(\overline{x(t)},\overline{y(t)})^\top$ is the analytic continuation
of $(x(t),y(t))^\top$ along
a loop contained in $D$, and the monodromy group of (9) is not
reducible for $\mu\not\in\Z$, then the solutions
$(\overline{x(t)},\overline{y(t)})^\top$ and $(x(t),y(t))^\top$
are linearly independent. This together with $y(t)\neq 0$ for $h\in
(1,\infty)$ (Proposition 2) shows that the imaginary part of
$F(t)$ has at most ${\rm deg}\,P$ zeros on $(1,\infty)$. Suppose finally
that $|t|$ is sufficiently big. As $|y(t)|\geq c|t|^{\lambda^*}$ then
$F(t) \lesssim  t^{s-\lambda^*}.$ Summing up the above information,
we obtain that the increase of the argument of $F(t)$, when running the
boundary of $D$ is at most $2\pi(1+ {\rm deg}\,P + s- \lambda^*)$.
Moreover, in the exceptional case when $\lambda^*=\frac12$, one has
$F(t)\sim ct^{s-\frac12}/\log t$ for large $|t|$, which yields
a stronger result: the increase of the argument of $F$ on $|t|=R$ is
strictly less than $2\pi(s-\frac12)$. Therefore the total increase of
the argument in this case is $<2\pi(1+ {\rm deg}\,P + s- \frac12)$.
This fact is useful only if $s-\frac12\in \N$ but we need it below.
One can deduce from the preceding proof of Proposition 4 that:
${\rm deg}\, P=[s-\lambda]$ if $\lambda-\mu\leq 1$,
${\rm deg}\, P=[s-\mu-2]$ if $\lambda-\mu>1$ and $s-\frac12$ are both
integers, and ${\rm deg}\, P=[s-\mu-1]$ otherwise.
(If ${\rm deg}\, P<0$, one takes $P\equiv 0$.) On its hand,
$\lambda^*=\mu$ if $\lambda-\mu\leq 1$ and $\lambda^*=\lambda-1$
otherwise. Therefore, by Proposition 4, the difference between the number
of zeros and poles in $D$ of the meromorphic function $F(t)=I(t)/y(t)$
is bounded by ${\rm dim}\, V_s-1$. By Proposition 3, this yields that
$I(t)$ has at most $[\lambda^*]+{\rm dim}\, V_s$ zeros in $D$.
Theorem 1 is proved. $\Box$

\vspace{3ex}
\noindent
{\large\bf 3. The applications}

\vspace{2ex}
\noindent
In this section we prove Theorems 2 and 3. Before that, let us point
out that some but not everything included in Table 1 is an evident fact.
However, since the procedure of deriving the related Fuchsian systems
is more or less known, we are not going to discuss in more details how
all these systems were obtained.

Given $i, j$ nonnegative integers,
denote $I_{ij}(h)=\int\!\!\int_{H<h}x^iy^jdxdy.$ Then
$${\cal V}_n=\{I(h)=\sum_{0\leq i+j\leq
n-1}c_{ij}I_{ij}(h)\}.\eqno(10)$$

\vspace{2ex}
\noindent
{\bf Lemma 1.} {\it Let ${\bf I}=(I_1,I_2)$ be as in Table 1.
Then for $n\geq 3$ one can express the function $I(h)$ from $(10)$ in
the form
$I(h)=\alpha(h)I_1(h)+\beta(h)I_2(h)$
where $\alpha(h)$ and $\beta(h)$ are polynomials of degrees as follows:}

\vspace{1ex}
(i) {\it ${\rm deg}\,\alpha=[\frac{n-1}{2}]$,
${\rm deg}\,\beta=[\frac{n-2}{2}]$ in cases $1)$ and} 2);

\vspace{1ex}
(ii) {\it ${\rm deg}\,\alpha=[\frac{n-1}{3}]$,
${\rm deg}\,\beta=[\frac{n-3}{3}]$ in case} 3);

\vspace{1ex}
(iii) {\it ${\rm deg}\,\alpha=[\frac{n-1}{2}]$,
${\rm deg}\,\beta=[\frac{n-3}{2}]$ in cases $4)$ and} 5);

\vspace{1ex}
(iv) {\it ${\rm deg}\,\alpha={\rm deg}\,\beta=[\frac{n-3}{2}]$
in cases $6)$ and} 7).

\vspace{1ex}
\noindent
{\it Moreover, the coefficients in $\alpha(h)$ and $\beta(h)$ may take
arbitrary values, except in case {\rm (iv)}. The dimension of the
vector space ${\cal V}_n$ in the case {\rm (iv)} equals
$[\frac{n-1}{2}]+[\frac{n-1}{4}] + 1$.}

\vspace{2ex}
\noindent
{\bf Proof.} For some of the cases, the results in Lemma 1 are already
known.
The result in case 1) was proved by Petrov in [5]. For cases 4) and 5)
see Petrov [8], [6], respectively.
The result for 2) follows from the considerations in [2] and [3].
The result concerning 3) is proved in [1].
Let us consider cases 6) and 7) from Table 1.
By symmetry, we have  $I_{ij}(h)=I_{ji}(h)$ and
$I_{ij}(h)\equiv 0$ whenever $i$ or $j$ is an odd number.
To establish the relations between the integrals $I_{ij}(h)$,
we take the equation $H\equiv x^2+y^2+\nu x^2y^2=h,\quad\nu=\pm1$ and
multiply both sides by the one-form $x^iy^{j+1}dx$. Afterwards
integrate the result along the oval $H=h$ and apply Green's
formula. One obtains the relation
$$(j+1)I_{i+2,j}+(j+3)I_{i,j+2}+\nu (j+3)I_{i+2,j+2}= (j+1)hI_{ij}.$$
Similarly, multiplying by $x^{i+1}y^jdy$ and integrating, we get
another relation
$$(i+3)I_{i+2,j}+(i+1)I_{i,j+2}+\nu (i+3)I_{i+2,j+2}= (i+1)hI_{ij}.$$
Combining these equations we easily obtain
$$\begin{array}{l}
\nu(i-j)I_{i+2,j+2}=(j+1)I_{i+2,j}-(i+1)I_{i,j+2},\quad i\neq j,\\[2mm]
\nu(i+3)I_{i+2,i+2}=-(2i+4)I_{i+2,i}+(i+1)hI_{ii},\\[2mm]
\nu(i+5)I_{i+4,0}=[\nu (i+2)h-1]I_{i+2,0}-3I_{i,2}+hI_{i,0}.
\end{array}\eqno(11)$$
For $i=j=0$, we get (noticing that $I_{00}=-I_1$ and
$I_{20}=I_{02}=-I_2$)
$$\begin{array}{l}
\textstyle I_{22}=\frac43\nu I_2-\frac13\nu hI_1\\[2mm]
\textstyle I_{40}=I_{04}=(-\frac25 h+\frac45\nu)I_2-\frac15\nu h I_1.
\end{array}$$
Then, using (11) with $i,j$ even, we easily prove the
assertion in (iv) by induction.

It remains to calculate the dimension of the vector space ${\cal V}_n$.
Clearly, we have ${\rm dim}\,{\cal V}_1=1$, ${\rm dim}\,{\cal V}_3=2$,
${\rm dim}\,{\cal V}_5=4$. By (11), the only new functions in
${\cal V}_{2m+1}$ (compared to ${\cal V}_{2m-1}$) are
$I_{2m,0}$ and, if $m$ is even,  $I_{m,m}$. Hence, the integrals
$I_{2k,0}$, $0\leq k\leq m$ and $I_{2k,2k}$, $1\leq k\leq m/2$ form
a basis in ${\cal V}_{2m+1}$.
These integrals are independent, since the leading term of
$I_{2k,0}$, $k\geq 2$ is proportional to $h^{k-1}(\nu I_1+2I_2)$
and the leading term of  $I_{2k,2k}$, $k\geq 1$ is proportional to
$h^kI_1$. The above argument implies that
${\rm dim}\,{\cal V}_n=[\frac{n-1}{2}]+[\frac{n-1}{4}] + 1$. $\Box$

\vspace{2ex}
\noindent
{\bf Remark 1.} In cases (i)-(iii) of Lemma 1, the result remains
true even for $n=1,2$, under the convention that a polynomial $\beta(h)$
of negative degree is taken to be zero. In case (iv), one has to take
$\beta(h)=0$, ${\rm deg}\,\alpha=0$ for $n=1,2$.

\vspace{2ex}
\noindent
{\bf Corollary 1.} {\it The dimension of the vector space ${\cal V}_n$,
$n\geq 1$, related to arbitrary polynomial perturbations of degree $n$
in $(3)$, in the cases $1)-7)$ of Table 1 is as follows:

\vspace{1ex}
$n,\qquad\qquad\,$ in cases $1)$ and $2)$;

\vspace{1ex}
$[\frac{2n+1}{3}],\qquad\;$ in case $3)$;

\vspace{1ex}
$2[\frac{n-1}{2}]+1,\;\;$ in cases $4)$ and $5)$;

\vspace{1ex}
$[\frac{n-1}{2}]+[\frac{n-1}{4}]+1,\;\;$ in cases $6)$ and $7)$;}

\vspace{2ex}
\noindent
{\bf Proof of Theorems 2 and 3.} Let us first note that $|\lambda-\mu|\leq
1$ for all cases 1)--7) in Table 1, which yields that $[\lambda^*]=0$.
We put

\vspace{1ex}
$s=\frac{n+1}{2}$ in cases 1) and 2),

\vspace{1ex}
$s=\frac{n}{3}+\frac12$ in case 3),

\vspace{1ex}
$s=[\frac{n+1}{2}]$ in cases  4) and 5),

\vspace{1ex}
$s=[\frac{n-1}{2}]+\frac12$ for $n\geq 3$, $s=1$ for $n=1,2$
in cases 6) and 7).

\vspace{1ex}
\noindent
It is easy to check that, with this choice of $s$,
${\cal V}_n\subset V_s$. For this purpose, one can perform the
inverse transformation $(I_1,I_2) \rightarrow (I_1, I_2-I_1/\omega)$
bringing $\bf A$ to a normal form (5) and then use the formulas
for the solution given in the proof of Proposition 4.
Hence, it suffices to verify in each case that ${\rm deg}\,\alpha
+\lambda\leq s$, ${\rm deg}\,\beta+\mu\leq s$
(the first inequality should be strong in cases 6) and 7);
see also Remark 1). Then we compare the dimensions of ${\cal V}_n$
and $ V_s$ (Proposition 4 and Corollary 1). One obtains that

$\dim {\cal V}_n = \dim  V_s $ in cases 1)--5), as well as in 6) and 7),
provided that $n\leq 6$,

$\dim V_s - \dim {\cal V}_n = [\frac{n-3}{4}]$ in
cases 6) and 7), if $n\geq 7$.

\noindent
Thus, the results follow from Theorem 1, taking into account that
$I(h)$ has always a zero at $h_0=0$. $\Box$

\vspace{2ex}
\noindent
{\bf Some other examples.}
Let us consider in brief system 8) from Table 1. Instead of
(3) and (4), we have
$$\begin{array}{l}
\dot{x}=H_y/M+\varepsilon f(x,y),\\
\dot{y}=-H_x/M+\varepsilon g(x,y),\end{array}\eqno(3')$$
where $H=x^{-3}(y^2-2x^2+x)$,
$M(x)=x^{-4}$, $f,g$ are polynomials of degree at most $n$, and
$$I(h)=\oint_{H=h}M(x)[g(x,y)dx-f(x,y)dy],\quad
h\in\Sigma=(-1,0).\eqno(4')$$
Note that in case 8), $\bf I$ is analytic in a neighbourhood of $h_0=-1$.
Define by ${\cal V}_n$ the linear space of integrals $(4')$ and let
${\cal D}=\C\setminus[0,\infty)$.

\vspace{2ex}
\noindent
{\bf Theorem 4.} {\it For  system $8)$ in Table 1,
the linear space of integrals ${\cal V}_n$ has a dimension $n+1$ and
is Chebyshev in ${\cal D}$, with accuracy as follows: one for $n=2$,
two for $n=1,3$, three for $n=0$ and $n-3$ for $n\geq 4$.}

\vspace{1ex}
\noindent
Taking $n=2$, we get the following result about the number of limit
cycles in $(3')$:

\vspace{1ex}
\noindent
{\bf Corollary 2.} {\it For any quadratic perturbation of the reversible
quadratic system $(3')$, the cyclicity of the period annulus around the center
at $(x,y)=(1,0)$ is two.}

\vspace{2ex}
\noindent
{\bf Proof of Theorem 4.} Denote $I_{kl}=\int\!\!\int_{H<h}M(x)x^ky^ldxdy$,
$I_k=\int_{H=h}M(x)x^{k-1}ydx$; thus ${\bf I}=(I_2,I_1)^\top$.
By symmetry, $I_{kl}=0$ for $l$ odd. In the same way as above,
we obtain the relations
$$\begin{array}{l}
I_{k,l+2}=\frac{2l+2}{2k+3l+3}(I_{k+2,l}-I_{k+1,l}),
\;k=-1, 0,\ldots,\;l=0, 2,\ldots,\\[2mm]
(k-\frac12)hI_{k+2}=(4-2k)I_{k+1}+(k-\frac72)I_k,\;k=0, 1, 2,\ldots,
\end{array}$$
and use the first of them to get the expression
$$I(h)=\sum_{k=0}^nc_kI_k(h),\quad c_k\; \mbox{\rm independent},
\eqno(12)$$
and then the second one to obtain
$$\begin{array}{ll}
I(h)=P_0I_1(h)+P_1(h)I_2(h), &{\rm for}\; n=0,1,2,\\[1mm]
I(h)=h^{-1}[P_1(h)I_1(h)+P_2(h)I_2(h)], &{\rm for}\; n=3,4,\\[1mm]
I(h)=h^{3-n}[P_{n-3}(h)I_1(h)+P_{n-2}(h)I_2(h)], &{\rm for}\; n\geq 5,
\end{array}\eqno(13)$$
where $P_k$ denotes a polynomial of degree $k$.
By (12), ${\rm dim}\,{\cal V}_n=n+1$. Given $n$, we choose
$s=\frac74$ if $n=0,1,2$, $s=\frac{11}{4}$ if $n=3,4$,
$s=n-\frac54$ if $n\geq 5$, and consider the corresponding linear
space $V_s$. Its dimension is ${\rm dim}\,V_{\frac74}=3$,
${\rm dim}\,V_{\frac{11}{4}}=5$,
${\rm dim}\,V_{n-\frac54}=2n-3$, respectively. For each $n$ the function
$I(h)$ in (13), multiplied by an appropriate power of $h$, belongs to the
respective $V_s$. The result then follows from Theorem 1. $\Box$

\vspace{1ex}
\noindent
Our last example is concerned with the Hamiltonian
$$H=x^2+y^2-x^4-ax^2y^2-y^4,\quad a>2,\eqno(14)$$
which comes from the cubic
Hamiltonian vector field having a rotational symmetry of order 4.
In complex coordinates $z=x+iy$, such a field is presented by
a complex equation
$\dot{z}=-iz+Az^2\overline{z}+B\overline{z}^3$, $A,B\in\C$, ${\rm Re}\,A=0.$
Take a polynomial perturbation in (3) which is semi-even with respect to $x$:
$$f(-x,y)=f(x,y), \quad g(-x,y)=-g(x,y),\quad {\rm deg}\, P,Q\leq n\eqno(15)$$
and consider the integral (4) where $\Sigma=(\frac{1}{a+2},\frac14)$ and
the integration is along the oval $\delta(h)\subset\{H=h\}$ surrounding
the center at $(\frac{1}{\sqrt2},0)$. As in Lemma 1, we can derive
relations between the integrals involved in (4) and then use them to
rewrite $I(h)$ in the form $I(h)=P(h)I_1(h)+Q(h)I_2(h)$
where  $I_1=\int_{\delta(h)}x^2dy$, $I_2=\int_{\delta(h)}x^2y^2dy$,
and $P,Q$ are polynomials with independent coefficients and degrees
$[\frac{n-2}{4}]$, $[\frac{n-4}{4}]$, respectively. The related
vector space ${\cal V}_n$ has a dimension $[\frac{n}{2}]$. The vector
function ${\bf I}=(I_1,I_2)^\top$ satisfies a system (2) with a matrix
(which is too large to fit in Table 1)
$${\bf A}=\left(\begin{array}{cc}
{\displaystyle\frac{4h-1}{3}} &
{\displaystyle\frac{a-2}{3}}\\[4mm]
{\displaystyle\frac{4h-1}{15(a+2)}} &
{\displaystyle\frac{4h}{5}+\frac{a-14}{15(a+2)}}
\end{array}\right)$$
Clearly, conditions (H1)--(H3) are satisfied with $h_0=\frac14$,
$h_1=\frac{1}{a+2}$. Denote ${\cal D}=\C\setminus(-\infty, h_1]$.
Take $s=\frac{n+1}{4}$, then evidently ${\cal V}_n=V_s$. Applying
Theorem 1, we obtain

\vspace{2ex}
\noindent
{\bf Theorem 5.} {\it For any system $(3)$ satisfying $(14)$ and
$(15)$, the linear space of integrals ${\cal V}_n$ has a dimension
$[\frac{n}{2}]$. Moreover, ${\cal V}_n$ is Chebyshev with accuracy
1 in ${\cal D}$ and it is Chebyshev in $\Sigma$.}

\vspace{2ex}
\noindent
Theorem 5 is useful for estimating the number of limit cycles not
surrounding the origin that are born in small semi-even polynomial
perturbations of the cubic Hamiltonian vector field with a rotational
symmetry of order 4.

{\large\bf  Appendix: Non-oscillation  and Sturm type theorems}

\vspace{2ex}
\noindent
The classical Sturm theorem can be used to find bounds for the number of
the zeros of the solutions
of linear non-autonomous differential equations  on a
real interval. In the context of the present paper a Sturm type
non-oscillation theorem was
recently proved by Petrov [8]. The proof uses of course  topological arguments.
It is natural
to ask whether the results of the present paper couldn't be deduced in such
a way. The answer turns
out to be negative in general, and our main Theorem 1 is essentially a
non-oscillation result in a
 complex domain. On the other hand our proofs also rely on topological
arguments: the argument
principle for real analytic functions in a complex domain. Therefore we may
call Theorem 1 a {\it Sturm type theorem in a complex domain}.

To compare these two approaches (real and complex) we
give below an example in which a  Sturm type theorem in a real domain
can still be proved. We shall follow closely Petrov [8]. As in the
introduction, $P$ and $Q$ are the polynomials from the definition of $V_s$
and it is assumed (for definiteness) that $h_0<h_1$.

\vspace{2ex}
\noindent
{\bf Theorem. } {\it Assume that conditions {\rm (H1)--(H3)} hold and
$2\lambda\not\in\Z$. Then any nontrivial function in $V_s$ has at most
${\rm deg}\,P+ {\rm deg}\,Q+1$ zeros in the interval  $(-\infty,h_0)$.
In particular, if $|\lambda-\mu|<1$, then $V_s$ is a Chebyshev vector
space in $(-\infty,h_0)$.}

\vspace{2ex}
\noindent
{\bf Proof.} As in section 2, it is sufficient to consider (2) as
a system for ${\bf I}(h)=(x(t),y(t))^\top$, with ${\bf A}$ taken in a
normal form (7). For $k$ a nonnegative integer, denote
$$\omega_{2k+1}=(k+\lambda)(k+\lambda-1),\;
\omega_{2k+2}=(k+\mu)(k+\mu-1),$$
$${\Omega}_k=\pmatrix{\omega_{2k+1} & 0\cr
0&\omega_{2k+2}},\qquad {\bf R}_k=k\pmatrix{
\lambda+k-1 & \mu\omega \cr \lambda/\omega & \mu+k-1}.$$
Following [8] we introduce the operator
$${\bf L}=\pmatrix{L&0\cr 0&L},\qquad L=t(t-1)\frac{d^2}{dt^2}.$$
By Proposition 1, we have $\quad {\bf L}{\bf I}=\Omega_0{\bf I}$.
Next, we prove that under hypotheses (H1)--(H2), the operator ${\bf L}$
satisfies also the following identities:
 $\;{\bf L}(t^k{\bf I})=t^k{\Omega}_k{\bf I}-t^{k-1}{\bf R}_k{\bf I}$,
$\; k\in\N$. Indeed, taking into account the form of the matrix in (7),
and denoting for short $\delta={\rm det}\,{\bf A}=t(t-1)/\lambda\mu$,
we obtain
$$\begin{array}{l}
{\bf L}(t^k{\bf I})=t(t-1)(t^k{\bf I})''\\
= t(t-1)[t^k{\bf I}''+2kt^{k-1}{\bf I}'+k(k-1)t^{k-2}{\bf I}]\\
= [t^k{\bf L}+2k\lambda\mu t^{k-1}\delta{\bf A}^{-1}
+k(k-1)(t^k-t^{k-1})]{\bf I}\\
= t^k[{\Omega}_0+2k\lambda\mu(\delta{\bf A}^{-1})'+k(k-1)]{\bf I}
+t^{k-1}[2k\lambda\mu(\delta{\bf A}^{-1})(0)-k(k-1)]{\bf I} \\
= t^k{\Omega}_k{\bf I}-t^{k-1}{\bf R}_k{\bf I}.\end{array}$$
Assume that $2\lambda$ is not integer. Then it is easy to
verify that
the constants $\omega_j$ are all different. This implies that there
exist scalar functions of the form
$$\begin{array}{l}
x_{k}(t)=[t^k+{\rm O}(t^{k-1})]x(t)+{\rm O}(t^{k-1})y(t),\\
y_{k}(t)={\rm O}(t^{k-1})x(t)+[t^k+{\rm O}(t^{k-1})]y(t),\end{array}$$
satisfying the equations
$$Lx_{k}(t)=\omega_{2k+1}x_{k}(t),\quad
Ly_{k}(t)=\omega_{2k+2}y_{k}(t),$$
where ${\rm O}(t^{k-1})$ denotes different polynomials of degree $k-1$.
To verify this, we ask for a ${\bf I}_k(t)=(x_k(t), y_k(t))^\top$
in the form
$${\bf I}_k=\sum_{j=0}^k {\bf B}_j t^j{\bf I},$$
${\bf B}_k$ the unity matrix, ${\bf B}_j$ to be determined for $j<k$.
As the operator  ${\bf L}$ commutes with the constant matrices, we have
$${\bf L}{\bf I}_k=\Omega_kt^k{\bf I}+ \sum_{j=0}^{k-1}
({\bf B}_j\Omega_j- {\bf B}_{j+1}{\bf R}_{j+1})t^j{\bf I}=
\Omega_k{\bf I}_k$$
and the matrices  ${\bf B}_j$, $j=k-1, k-2,\ldots, 0$ are determined
recursively from the equations
$${\bf B}_j\Omega_j-\Omega_k{\bf B}_j={\bf B}_{j+1}{\bf R}_{j+1}$$
which is possible because $\omega_j$ are all different.
Therefore, there is a basis in the space $V_s$ consisting
of the eigenfunctions of the operator $L$. Taking into account that,
by Proposition 2, $x_k$ and $y_k$ do not vanish for $t<0$, we apply
Petrov's elimination technique [8] to prove that any function
in $V_s$ has at most ${\rm deg}\,P+ {\rm deg}\,Q+1$ isolated zeros in
$(-\infty,0)$.  Especially, in the case when $|\lambda-\mu|<1$,
this means that $V_s$ is Chebyshev in $(-\infty,0)$. $\Box$

\vspace{1ex}
Note that the above proof works only on the open intervals
having $h_0$ as an endpoint and where ${\rm det}\,{\bf A}$ is
positive (because the Sturm theorem applies in a
backward direction here). Also note that when $|\lambda-\mu|>1$,
the above estimate, although it concerns the interval
$(-\infty,h_0)$ only, is weaker than the estimate for the
whole $\cal D$ obtained in Theorem 1.

\vspace{4ex}
\noindent
{\bf References}

\small
\newcounter{num}
\begin{list}{[{\arabic{num}}]}
{\usecounter{num}\setlength{\itemsep}{2.0mm}\setlength{\parsep}{0.0mm}
\setlength{\itemindent}{-2.0mm}}
\item {L. Gavrilov}, Nonoscillation of elliptic integrals related to
cubic
polynomials with symmetry of order three, {\em Bull. Lond. Math. Soc.}
{\bf 30} (1998), 267--273.
\item {L. Gavrilov}, Abelian integrals related to Morse polynomials and
perturbations of plane Hamiltonian vector fields,
{\em Ann. Inst. Fourier, Grenoble}
{\bf 49} (1999), 611--652.
\item Emil Horozov, Iliya D. Iliev, Linear estimate for the number of
zeros of
Abelian integrals with cubic Hamiltonians, {\em Nonlinearity} {\bf 11}
(1998),
no. 6, 1521--1537.
\item K. Iwasaki, H. Kimura, Sh. Shimomura, M. Yoshida, From Gauss to
Painlev\'e. A modern theory of special functions.
 Aspects of Mathematics, vol. 16, Vieweg, Braunschweig (1991).
\item {G.S. Petrov},  Elliptic integrals and their non-oscillation,
{\em Funct. Anal. Appl.} {\bf 20} (1986), no 1, 46--49.
\item {G.S. Petrov},    Complex zeroes of an elliptic integral,
{\em Funct. Anal. Appl.} {\bf 21} (1987), no 3, 87--88.
\item {G.S. Petrov},    Nonoscillation of elliptic integrals,
{\em Funct. Anal. Appl.} {\bf 24} (1990), no 3, 45--50.
\item {G.S. Petrov},  On the nonoscillation of elliptic integrals,
{\em Funct. Anal. Appl.} {\bf 31} (1997), no 4, 47--51.
\item {G.S. Petrov},  Complex zeroes of an elliptic integral,
{\em Funct. Anal. Appl.} {\bf 23} (1989), no 2, 88--89.
\item {G.S. Petrov},  The Chebyshev property of elliptic integrals,
{\em Funct. Anal. Appl.} {\bf 22} (1988), no 1, 83--84.
\item {C. Rousseau} and {H. \.{Z}o{\l}adek}, Zeros of complete elliptic
integrals for 1:2 resonance, {\em J. Differential Equations} {\bf 94}
(1991),
no 1, 41--54.
%
%
\item {\it Encyclopaedia of Mathematics}, vol. 7, pp. 30--33;
vol. 9, pp. 302--303, Kluwer (1991).
\end{list}
\end{document}